\documentclass{amsart}
\usepackage{amscd, graphics}

\newtheorem{thm}{Theorem}[section]

\newtheorem{prop}[thm]{Proposition}

\theoremstyle{definition}
\newtheorem{defn}[thm]{Definition}

\theoremstyle{remark}

\newtheorem{rem}{Remark}

\numberwithin{equation}{section}

\newcommand{\C}{\mathbb{F}}

\newcommand{\mcA}{\mathcal{A}}
\newcommand{\mcR}{\mathcal{R}}
\newcommand{\mcF}{\mathcal{F}}

\newcommand{\autcn}{{\rm Aut}({\mathbb C}^{n})}

\newcommand{\Cn}{{\mathbb C}^{n}}

\newcommand{\End}{\text{End }}

\begin{document}

\title[Degenerations of $R$-operators]{Degenerations and 
representations of twisted Shibukawa-Ueno $R$-operators}
\author{Robin Endelman}
\address{University of California at Davis, Davis, CA 95616,
U.S.A.}
\email{endelman@math.ucdavis.edu}
\author{Timothy J. Hodges}
\address{University of Cincinnati, Cincinnati, OH 45221-0025}
\email{timothy.hodges@uc.edu}
\thanks{The first author was supported in part by NSF VIGRE Grant No. DMS-0135345.
The second author was supported in part by NSA grant 
MDA904-99-1-0026 and by the Charles P. Taft Foundation}

\begin{abstract}
	We study degenerations of the Belavin $R$-matrices via the infinite dimensional
 operators defined by Shibukawa-Ueno. We define a two-parameter family of 
	generalizations of the Shibukawa-Ueno $R$-operators.  These 
	operators have finite dimensional representations which 
	include Belavin's $R$-matrices in the elliptic case, 
	a two-parameter family of twisted affinized Cremmer-Gervais 
	$R$-matrices in the trigonometric case, and a two-parameter 
	family of twisted (affinized) generalized Jordanian 
	$R$-matrices in the rational case.  We find finite dimensional 
	representations which are compatible with the elliptic to 
	trigonometric and rational degeneration.  We further show that certain
	members of the elliptic family of operators have no finite dimensional 
	representations.  These $R$-operators 
	unify and generalize earlier constructions of Felder and 
	Pasquier, Ding and Hodges, and the authors, and illuminate 
	the extent to which the Cremmer-Gervais $R$-matrices (and their rational 
	forms) are degenerations of Belavin's $R$-matrix.

\end{abstract}
\maketitle

\section{Introduction}
	In this article, we give a unified description of finite 
	dimensional representations of twisted Shibukawa-Ueno 
	$R$-operators.  In \cite{SU}, Shibukawa and Ueno described a set 
	of solutions of the Yang-Baxter equation on the field of 
	meromorphic functions on $\C$.  These solutions were of three 
	types, elliptic, trigonometric and rational.  In \cite{FP}, Felder 
	and Pasquier showed that the elliptic solutions could be twisted 
	and restricted to a finite dimensional subspace in such a way that 
	they produced Belavin's solutions of the Yang-Baxter equation.  
	Ding and the second author observed in \cite{DH} that a similar 
	procedure applied to the analogous trigonometric solutions of the 
	constant Yang-Baxter equation yielded the Cremmer-Gervais 
	$R$-matrices.  The authors then showed in \cite{EH1} that this 
	procedure applied to the constant rational solutions yielded a 
	generalization of the Jordanian $R$-matrix which quantized 
	certain solutions of the classical Yang-Baxter equation studied by 
	Gerstenhaber and Giaquinto in \cite{GG1,GG2}.  Our aim here is to 
	unify and generalize these constructions in order to explain to 
	what extent the Belavin $R$-matrices degenerate into affinized 
	Cremmer-Gervais $R$-matrices and their rational analogs.

First we observe that a simple twisting procedure enables us to 
extend the Shibukawa-Ueno operators to a two parameter family of 
Yang-Baxter operators. We then look for finite dimensional 
representations of these operators which are compatible with their 
degeneration. In the elliptic case, only certain members of this 
family have such finite dimensional representations and these 
representations are essentially those studied by Felder and Pasquier. 
In the trigonometric case the whole family has a simultaneous finite 
dimensional representation and the family of $R$-matrices obtained is 
a 2-parameter family of deformations of the affinized Cremmer-Gervais 
$R$-matrices constructed using an analogous twisting mechanism. The 
Belavin $R$-matrices then degenerate naturally not to the 
Cremmer-Gervais matrices themselves but to other members of this 
family. This is analogous to the fashion in which Antonov, Hasegawa and 
Zabrodin realised the Cremmer-Gervais matrices as twists of 
degenerations of the Belavin matrices \cite{AHZ}. In the rational 
case, the generalized Shibukawa-Ueno operators also have a 
simultaneous finite dimensional representation which yields a two 
parameter family of deformations of the (affinized) generalized Jordanian $R$-matrices 
constructed by the authors in \cite{EH1}.

In order to obtain the required degenerations of the 
finite dimensional matrices, we degenerate the finite dimensional 
subspaces in a fashion consistent with the degeneration of the 
operators.  Thus the subspace on which we represent the elliptic 
operator has a basis of elliptic functions, while the subspace on 
which we represent the trigonometric operators and rational operators 
have bases of trigonometric and rational functions, respectively. 

The authors would like to thank Jintai Ding for many helpful discussions.  

\section{Background}
\subsection{Yang-Baxter equation on fields of meromorphic functions}

We begin by setting out an appropriate context in which to discuss the 
Yang-Baxter equation on function fields. Denote the field of 
meromorphic functions on ${\mathbb C}^{n}$ by ${\mathcal M}({\mathbb 
C}^{n})$. In general it is not possible to extend a linear operator 
$\mcR \in \End ({\mathcal M}({\mathbb C}^{2}))$ to an operator $\mcR_{12} 
\in \End ({\mathcal M}({\mathbb C}^{3}))$. Thus we must work inside a 
subring of $\End ({\mathcal M}({\mathbb C}^{2}))$ in which this is 
possible.

Denote by $\autcn$ the group of automorphisms of $\Cn$. 
For any function  $\phi \in \autcn $ define
$\phi^{*} \in \End( {\mathcal M}({\mathbb C}^{n}))$ by 
$\phi^{*}(f) = f \circ \phi$. Set $G(n) = \{\phi^{*}| \phi \in 
{\rm Aut}({\mathbb C}^{n})\}$, and let 
${\mathcal A(n)}$ denote the subalgebra 
of ${\rm End}({\mathcal M}({\mathbb C}^{n}))$ 
generated by $G(n)$ over the subfield ${\mathcal M}({\mathbb C}^{n})$ 
(acting as multiplication operators).  It is easily verified that 
$G(n)$ is linearly independent over ${\mathcal M}({\mathbb C}^{n})$, 
and that ${\mathcal M}(\Cn)$ is a $G(n)$-module algebra.  For 
$\phi^{*} \in G(n)$ and $f \in {\mathcal M}(\Cn)$, 
$\phi^{*} \cdot f = \phi^{*}(f) \phi^{*}$, and so 
${\mathcal A(n)}$ is isomorphic to the smash product 
${\mathcal M}({\mathbb C}^{n}) \# G(n)$. Denote ${\mathcal A(2)}$ by 
${\mathcal A}$.

	For $\phi \in {\rm Aut}({\mathbb C}^{2})$, define $\phi_{ij} \in 
{\rm Aut}({\mathbb C}^{3})$ by $\phi$ acting on the $i^{\rm th}$ and $j^{\rm th}$ 
variables.  Then $(\phi^{*})_{ij} = (\phi_{ij})^{*} \in G(3)$, and for 
$\mcR = \sum_{\alpha} f_{\alpha}(z_{1},z_{2})\phi_{\alpha}^{*}\in 
{\mathcal A}$, with 
$f_{\alpha} \in 
{\mathcal M}({\mathbb C}^{2})$ and $\phi_{\alpha}^{*} \in G(2)$, we may define 
$\mcR_{ij} \in {\mathcal A}(3)$ by 
$$\mcR_{ij} = \sum_{\alpha} 
f_{\alpha}(z_{i},z_{j})(\phi_{\alpha}^{*})_{ij}.$$

\begin{defn}  
A solution, $\mcR$, in ${\mathcal A}$ to the Yang-Baxter equation will be called an 
{\em $R$-operator}.  
As in the case for $R$-matrices, spectral parameter-dependent 
$R$-operators are maps $\mcR : \Omega  
\longrightarrow {\mathcal A}$, where $\Omega \subset {\mathbb C}^{k}$ 
with $k = 
1$ or $2$, satisfying the spectral parameter-dependent Yang-Baxter equation 
$$\mcR_{12}(\lambda_{1},\lambda_{2})\mcR_{13}(\lambda_{1},\lambda_{3}) 
\mcR_{23}(\lambda_{2},\lambda_{3}) = \mcR_{23}(\lambda_{2},\lambda_{3}) 
\mcR_{13}(\lambda_{1},\lambda_{3}) \mcR_{12}(\lambda_{1},\lambda_{2}).$$
\end{defn}

Denote by $P$ the $R$-operator, in $G(2)$, which acts on 
${\mathcal M}({\mathbb C}^{2})$ by $P \cdot f(z_{1},z_{2}) = f(z_{2},z_{1})$.  

\subsection{The Shibukawa-Ueno $R$-operators}

In \cite{SU}, Shibukawa and Ueno constructed a family of operators 
$\mcR^{\theta}(\lambda) \in {\mathcal A}$, depending on a holomorphic 
function 
$\theta$:
$$\mcR^{\theta}(\lambda) = G_{\theta}(z_{1}-z_{2},\lambda)P - 
G_{\theta}(z_{1}-z_{2},\kappa)I,$$
where $\displaystyle{G_{\theta}(z,\lambda) = \frac{\theta'(0) 
\theta(z+\lambda)}
{\theta(z) \theta(\lambda)}}$ and 
$I$ is the identity operator. These operators 
satisfy the Yang-Baxter equation for any scalar $\kappa$ and any 
function $\theta$ satisfying the ``three-term equation":
$$\sum_{\stackrel{{\rm cyc~perms}}{{\rm of~} y,z,w}}
\theta(x+y)\theta(x-y)\theta(z+w)\theta(z-w) = 0.$$
Analytic solutions of the three-term equation 
are \cite{WW} $\theta(z) = Ae^{Bz^{2}}\vartheta_{1}(Cz,\tau)$ 
and its trigonometric and rational degenerations, 
where $A$, $B$, $C$ are constants, 
$\tau \in {\mathbb H} = \{z \in {\mathbb C} ~|~ 
{\rm Im} z > 0\}$, and 
$\vartheta_{1}$ is Jacobi's first theta function:
$$\vartheta_{1}(z,\tau) = -\sum_{m \in {\mathbb Z}}{\rm exp}\left\{ 
	\pi i \left(m + \frac{1}{2}\right)^{2}\tau + 2 \pi i \left(m + 
	\frac{1}{2}\right)\left(z + \frac{1}{2}\right)\right\}.$$

	Here we will consider the particular cases when $A=1$ and 
$B=0$. Specifically, we consider the following families of 
$R$-operators:

\begin{itemize}
	\item[] Elliptic case: $\mcR^{e}(\lambda)$ with 
	$\theta(z) = \vartheta_{1}(z,\tau)$, $\tau \in {\mathbb H}$
	\item[] Trigonometric case: $\mcR^{t}(\lambda)$ with 
	$\theta(z) = \sin \pi z/\tau_{1}$, $\tau_{1} \in {\mathbb C} 
	\setminus \{0\}$ 
	\item[] Rational case: $\mcR^{r}(\lambda)$ with $\theta(z) = z$.
\end{itemize}

In \cite{SU}, Shibukawa and Ueno also introduce the notion of 
obtaining $R$-matrices by restriction of $\mcR^{\theta}(\lambda)$ 
to finite dimensional invariant subspaces.  
For $n < \infty$, let $V = \bigoplus_{a} {\mathbb C}f_{a}$ be 
an $n$-dimensional subspace of ${\mathcal M}({\mathbb C})$, with 
basis $\{f_{a} \mid a = 0, 1, \ldots, n-1\}$, and 
identify $V \otimes V$ with the space of functions in 2 variables 
${\rm Span}_{\mathbb C} 
\{f_{a}(z_{1})f_{b}(z_{2})\}$.  

\begin{defn} \cite{SU} If $\mcR \in {\mathcal A}$ is a solution of the 
Yang-Baxter equation and 
$V \otimes V$ is invariant under $\mcR$, then we say that 
$R|_{V \otimes V} \in {\rm End}(V \otimes V)$ is 
a {\em finite dimensional representation of $\mcR$}. 
 For a spectral parameter-dependent $R$-operator, $\mcR: \Omega 
\longrightarrow {\mathcal A}$, we say that $R|_{V \otimes V}$  is a finite 
dimensional representation of $\mcR$ if $V \otimes V$ is invariant under 
$\mcR(\lambda)$ for all $\lambda \in \Omega$.
\end{defn}

If $\mcR$ is an $R$-operator, then $\mcR|_{V\otimes V}$ is a 
matrix solution of the Yang-Baxter equation.  In Section 4, we 
give finite dimensional representations for twisted 
Shibukawa-Ueno $R$-operators defined in Section 3, in each of the 
elliptic, 
trigonometric, and rational cases, and the corresponding degenerations. 

\section{Twisted Shibukawa-Ueno operators}

The twists that we shall be considering are all of one simple kind.

\begin{thm}
Let $R(\lambda) \in {\rm End}(V \otimes V)$ be a solution of the Yang-Baxter equation. 
Let $B \colon 
(\mathbb{C},+) \to GL(V)$ be a homomorphism such that $R(\lambda)$ 
commutes with $B(\mu)\otimes B(\mu)$ for all $\lambda$ and $\mu$.  
For $\alpha,\beta \in {\mathbb C}$, set
$$
	F_{\alpha,\beta}(\lambda) = B(\alpha \lambda - \beta) \otimes 
	B^{-1}(\alpha \lambda - \beta)
$$
Then
$$
	R_{\alpha,\beta}(\lambda) = F_{\alpha, \beta}(-\lambda) 
	R(\lambda) F_{\alpha, \beta}(\lambda)
$$
also satisfies the Yang-Baxter equation.
\end{thm}

\begin{proof}
It suffices to prove the result in two separate case: 1) 
when $\alpha = 1$ and $\beta=0$ and 2) when $\alpha = 0$.
Let $B_{i}$ denote $B$ acting in the $i^{\rm th}$ component of the 
tensor product.  

Case 1: Suppose $\alpha = 1$ and $\beta=0$. Since 
$R_{1,0}(\lambda) = B_1(-2\lambda)R(\lambda)B_1(2\lambda) $, it 
actually suffices to prove that 
$\tilde{R}(\lambda) = B_1(\lambda)R(\lambda)B_1(\lambda)$ satisfies the 
Yang-Baxter equation. 

Now
\begin{alignat*}{2}
&&&	\tilde{R}_{12}(\lambda-\lambda')\tilde{R}_{13}(\lambda)\tilde{R}_{23}(\lambda')\\
&= &&B_1(\lambda-\lambda')R_{12}(\lambda-\lambda')B_1(\lambda'-\lambda)
B_1(\lambda)R_{13}(\lambda)B_1(-\lambda)B_2(\lambda')R_{23}(\lambda')B_2(-\lambda')\\
&= &&B_1(\lambda-\lambda')R_{12}(\lambda-\lambda')B_1(\lambda')
B_2(\lambda')R_{13}(\lambda)R_{23}(\lambda')B_2(-\lambda')B_1(-\lambda)\\
&= &&B_1(\lambda)B_2(\lambda')R_{12}(\lambda-\lambda')R_{13}(\lambda)R_{23}(\lambda')B_2(-\lambda')B_1(-\lambda)\\
&= &&B_1(\lambda)B_2(\lambda')R_{23}(\lambda')R_{13}(\lambda)R_{12}(\lambda-\lambda')B_2(-\lambda')B_1(-\lambda)\\
&= &&B_1(\lambda)B_2(\lambda')R_{23}(\lambda')R_{13}(\lambda)B_1(-\lambda')B_2(-\lambda')R_{12}(\lambda-\lambda')B_1(\lambda'-\lambda)\\
&=
&&B_2(\lambda')R_{23}(\lambda')B_2(-\lambda')B_1(\lambda)R_{13}(\lambda)B_1(-\lambda)B_1(\lambda-\lambda')
R_{12}(\lambda-\lambda')B_1(\lambda'-\lambda)\\
&= &&\tilde{R}_{23}(\lambda')\tilde{R}_{13}(\lambda)\tilde{R}_{12}(\lambda-\lambda')
\end{alignat*}

Case 2: Suppose $\alpha = 0$. Since $R_{0,\beta}(\lambda) = B_1(-2\beta)R(\lambda) 
B_2(2\beta) $, it suffices to 
prove that $\tilde{R}(\lambda) = B_1(\beta)) R(\lambda) (B_2(-\beta) $ satisfies the 
Yang-Baxter equation. 
Now
\begin{alignat*}{2}
&&&	\tilde{R}_{12}(\lambda-\lambda')\tilde{R}_{13}(\lambda)\tilde{R}_{23}(\lambda')\\
&= &&B_1(\beta)R_{12}(\lambda-\lambda')B_2(-\beta)
B_1(\beta)R_{13}(\lambda)B_3(-\beta)B_2(\beta)R_{23}(\lambda')B_3(-\beta)\\
&= &&B_1(\beta)R_{12}(\lambda-\lambda')B_2(\beta)
B_1(\beta)R_{13}(\lambda)B_3(-\beta)B_2(-\beta)R_{23}(\lambda')B_3(-\beta)\\
&= &&B_1(2\beta)B_2(\beta)R_{12}(\lambda-\lambda')R_{13}(\lambda)R_{23}(\lambda')B_2(-\beta)B_3(-2\beta)\\
&= &&B_1(2\beta)B_2(\beta)R_{23}(\lambda')R_{13}(\lambda)R_{12}(\lambda-\lambda')B_2(-\beta)B_3(-2\beta)\\
&= &&B_2(\beta)R_{23}(\lambda')B_1(2\beta)R_{13}(\lambda)B_3(-2\beta)R_{12}(\lambda-\lambda')B_2(-\beta)\\
&= &&B_2(\beta)R_{23}(\lambda')B_3(-\beta)B_1(\beta)R_{13}(\lambda)B_3(-\beta) B_1(\beta)R_{12}(\lambda-\lambda')B_2(-\beta)\\
&= &&\tilde{R}_{23}(\lambda')\tilde{R}_{13}(\lambda)\tilde{R}_{12}(\lambda-\lambda')
\end{alignat*}

\end{proof}

Suppose that $V$ has a basis $e_i$, for $i=1, \dots, n$. An 
operator $R(\lambda) \in {\rm End}(V \otimes V)$ is said to be 
homogeneous if it is of the form
$$
R(\lambda)(e_i \otimes e_j) = \sum_{k} a_k(\lambda) e_k \otimes 
e_{i+j-k}
$$

The diagonal operator $B(\mu)e_k = \exp(c\mu k) e_k$ satisfies the 
hypothesis of the theorem for any homogeneous operator $R(\lambda)$ 
and the resulting twisted operator is

$$
R_{\alpha, \beta}(\lambda)(e_i \otimes e_j) = \sum_{k} 
\exp\{2c[\alpha \lambda (i-k) - \beta (k-j)]\} 
a_k(\lambda) e_k \otimes e_{i+j-k}.
$$

The Cremmer-Gervais $R$-matrices \cite{CG} are homogeneous, given by
$$
	(R_{CG})_{ij}^{kl} = p^{2(j-k)}
		\begin{cases}
			q &  i=k \geq j=l,  \\
			q^{-1} & i=k < j=l, \\
			-\widehat{q} & i<k<j, i+j=k+l,  \\
			\widehat{q} & j \leq k < i, i+j=k+l, \\
			0 & {\rm otherwise},
		\end{cases}.
		$$
where 
$p^{n} = q $ and $\hat{q} = q - q^{-1}$.  
The 2-parameter versions ($p$ and $q$ independent) described in \cite{Hcg} are 
constructed from twists of this form ($\alpha = 0$).  

The corresponding solution 
of the braid equation $\check{R}_{CG}$ satisfies the Hecke relation 
$(\check{R}_{CG} -q)(\check{R}_{CG} + q^{-1}) = 0$ 
and as such, $R_{CG}$ has a standard affinization \cite[page 296]{KS} given by 
$R_{CG}(\lambda) = \hat{q} \eta P - \hat{\eta} R_{CG}$, 
where $\eta = e^{\pi i \lambda}$, whose matrix coefficients 
are given by
$$
	R_{CG}(\lambda)_{ij}^{kl} = p^{2(j-k)}
	\begin{cases}
		\widehat{q\eta^{-1}} &  i=j=k=l,  \\
		-\hat{\eta}q^{\text{sgn}(i-j)} & i=k \neq l=j, \\
		\text{sgn}(j-i) \hat{\eta}\hat{q}  & \min(i,j) < k < \max(i,j), i+j=k+l,  \\
		\hat{q} \eta^{\text{sgn}(j-i)} & j=k \neq l=i, \\
		0 & \text{otherwise}.
		\end{cases}
		$$

With $B(\mu) e_k = e^{2\pi i \mu k} e_k$, we define a two-parameter family of 
twisted Cremmer-Gervais operators by

\begin{equation}
R_{CG(\alpha,\beta)}(\lambda) = F_{\alpha,\beta}(-\lambda)R_{CG}(\lambda)
F_{\alpha,\beta}(\lambda).
	\label{eq:CGab}
\end{equation}
The matrix coefficients are then given by
$$ R_{CG(\alpha,\beta)}(\lambda)_{ij}^{kl} =
\zeta^{2(i-k)}\gamma^{2(j-k)}R_{CG}(\lambda)_{ij}^{kl}
$$
where $\zeta = e^{2 \pi i \alpha \lambda}$ and $\gamma = e^{2 \pi i \beta}$.  

An analogous version of the twisting theorem for operators on function spaces 
is the following.

\begin{thm}
Let $\mcR(\lambda) \in \mcA$ be an $\mcR$-operator. 
Let $\phi: {\mathbb C} \rightarrow {\rm Aut}({\mathbb C})$ be 
such that $\mcR(\lambda)$ 
commutes with $\phi^*(\mu)\otimes \phi^{*}(\mu)$ for all $\lambda$ and 
$\mu$.  For $\alpha, \beta \in {\mathbb C}$, set
$$
	\mcF_{\alpha,\beta}(\lambda) = \phi^*(\alpha \lambda - \beta) \otimes 
	(\phi^*)^{-1}(\alpha \lambda - \beta)
$$
Then
$$
	\mcR_{\alpha,\beta}(\lambda) = \mcF_{\alpha, \beta}(-\lambda) 
	\mcR(\lambda) \mcF_{\alpha, \beta}(\lambda)
$$
also satisfies the Yang-Baxter equation.
\end{thm}

\begin{proof} Analogous to the proof in the finite dimensional case.
\end{proof}

Define $\phi \colon (\mathbb{C}, +) \to 
{\rm Aut}(\mathbb{C})$ 
by $\phi(\lambda) \cdot z = z + \lambda$. Define 
also operators  $\tilde{\mcF}_{s} \in \mcA$ by $\tilde{\mcF}_{s} \cdot 
f(z_{1},z_{2}) = 
f(z_{1} + s, z_{2} - s)$. 
If $\mcR^{\theta}(\lambda)$ is the Shibukawa-Ueno $R$-operator defined 
above, then it is easily seen that $\mcR^{\theta}(\lambda)$ commutes 
with $\phi^{*}(\mu) \otimes \phi^{*}(\mu)$ and 
hence, 
\begin{eqnarray*}
	\mcR^{\theta}_{\alpha,\beta}(\lambda) & = & 
\mcF_{\alpha,\beta}(-\lambda)\mcR^{\theta}(\lambda)\mcF_{\alpha,\beta}(\lambda)\\
& = & G(z_{1} - z_{2} - 2(\alpha \lambda + \beta),\lambda) 
	 \tilde{\mcF}_{- 2\alpha \lambda}P 
	 - G(z_{1} - z_{2} - 2(\alpha \lambda + \beta),\kappa) 
	 \tilde{\mcF}_{-2\beta} 
\end{eqnarray*}
satisfies the Yang-Baxter equation.

\begin{rem}
\par Taking $\alpha = 1/2n$ and $\beta = \kappa/2n$ 
		($n \in {\mathbb Z}^{\times}$), 
		the resulting twisted $R$-operator 
		in the elliptic case is the same as that of 
		Felder and Pasquier (in \cite{FP}, they 
		twist the elliptic Shibukawa-Ueno $R$-operator by a 
		two-parameter $F(\lambda_1,\lambda_2)$).  
\par Taking $\alpha = 0$ and $\beta = \kappa/2n$, 
		$\mcF_{0,\kappa/2n}$ is the constant 
		twist applied to $\mcR^{t}_{\infty} 
		( = \lim_{{\rm Im}~\lambda \to \infty} \mcR^{t}(\lambda))$ 
		in \cite{DH} in the trigonometric case, 
		and to $\mcR^{r}_{\infty} ( = \lim_{\lambda \to 
		\infty} \mcR^{r}(\lambda))$ in \cite{EH1} in the rational case. 
\end{rem}

\section{Finite dimensional representations and Degenerations}

\subsection{Elliptic case}
	
Let ${\mathcal H}({\mathbb C}^{n})$ denote the space of holomorphic functions 
on ${\mathbb C}^{n}$.  
Fix $\tau \in {\mathbb H}$, the upper half-plane, and let 
\begin{equation} 
V^{e}_{n} \stackrel{\rm def}{=} \{f \in {\mathcal H}({\mathbb C}) ~|~ 
f(z+1) = (-1)^{n-1}f(z), f(z + \tau) = e^{-\pi i n \tau - 2\pi i n 
z}f(z)\}.
	\label{eq:Ve}
\end{equation}

\begin{prop}\label{ebasis} The vector space $V^{e}_{n}$ has 
basis ${\mathcal B} = \{\psi_{a}(z) | a \in {\mathbb Z}_{n}\}$ where
$$\psi_{a}(z) = \sum_{\stackrel{m \in {\mathbb Z}}{m \equiv a ({\rm 
mod}~n)}} e^{\pi i (m - \frac{n-1}{2})^{2}\tau/n + 2\pi i 
(m-\frac{n-1}{2}) z}.$$
\end{prop}

\begin{proof}  It is well known that ${\rm dim}(V^{e}_{n}) = n$, so it 
suffices to show that ${\mathcal B}$ spans $V^{e}_{n}$.  
%Note that ${\cal B}$ is well-defined:  $\psi_{a+n} = \psi_{a}$.  
Let $f \in V^{e}_{n}$, and set $g(z) = e^{\pi i (n-1)z}f(z)$.  
Since $f(z+1) = (-1)^{n-1} f(z)$, $g$ satisfies $g(z+1) = g(z)$, 
and has Laurent expansion 
$g(z) = \sum_{m \in {\mathbb Z}} a_{m}e^{2\pi i m z}.$
Thus 
$$f(z) = e^{-\pi i (n-1)z}g(z) = 
\sum_{m \in {\mathbb Z}} a_{m}e^{2\pi i (m - \frac{n-1}{2}) z}.$$
Now, writing $a_{m} = 
c_{m}e^{\pi i (m - \frac{n-1}{2})^{2}\tau/n}$, from the 
quasi-periodicity of $f$ with respect to $\tau$, 
we obtain $c_{m} = c_{m+n}$ for 
all $m \in {\mathbb Z}$.  Hence, 
%$f$ is determined by $c_{0}, \ldots, c_{n-1}$:
$f(z) = \sum_{a \in {\mathbb Z}_{n}} c_{a}\psi_{a}(z)$.
\end{proof}

Identify $V^{e}_{n} \otimes V^{e}_{n}$ with the function space 
${\rm Span}_{\mathbb C} 
\{\psi_{a}(z_{1})\psi_{b}(z_{2})\} 
= \{f \in {\mathcal H}({\mathbb C}^{2}) \mid f(-,z_{2}),f(z_{1},-) \in 
V^{e}_{n}\}$. Let $\theta_{a,b}$ denote the standard theta 
functions of rational characteristic, defined, for $a,b \in {\mathbb 
Q}$, by: 
$$\theta_{a,b}(z,\tau) = 
\sum_{m \in {\mathbb Z}} {\rm exp}\left\{\pi i \left(m + 
a\right)^{2}\tau + 2\pi i (m + a)(z + b)\right\}.$$

\begin{thm}  With $V^{e}_{n}$ defined as above (\ref{eq:Ve}), 
$V^{e}_{n} \otimes V^{e}_{n}$ is a finite dimensional 
representation of $\mcR^{e}_{1/2n,\kappa/2n}(\lambda) \in {\mathcal A}$, on which 
$\mcR^{e}_{1/2n,\kappa/2n}(\lambda)$ is equivalent to a Belavin 
$R$-matrix.  The matrix coefficients of 
$\mcR^{e}_{1/2n,\kappa/2n}(\lambda)|_{V^{e}_{n} \otimes V^{e}_{n}}$ with respect 
to the 
basis $\{\psi_{a} \otimes \psi_{b}\}$ are given by
$$R(\lambda)_{ij}^{kl} = \delta_{i+j,k+l (n)} ~
\frac{\theta_{\frac{1}{2}, \frac{1}{2}}'(0,n\tau) 
\theta_{\frac{i-j}{n}+ \frac{1}{2},\frac{1}{2}}(\lambda - 
\kappa,n\tau)}{\theta_{\frac{i-k}{n}+ 
\frac{1}{2},\frac{1}{2}}(-\kappa,n\tau) \theta_{\frac{k-j}{n}+ 
\frac{1}{2}, \frac{1}{2}}(\lambda,n\tau)}.$$
\end{thm}

\begin{proof}
	Set $\theta(z) = \vartheta_{1}(z,\tau) = 
-\theta_{\frac{1}{2},\frac{1}{2}}(z,\tau)$, and 
$\mcR(\lambda) = \theta(\lambda)\mcR^{e}_{1/2n,\kappa/2n}(\lambda)$:

\begin{eqnarray*}
	\mcR(\lambda) & = & \theta(\lambda)\mcF_{1/2n,\kappa/2n}(-\lambda)
	\mcR^{e}(\lambda) 
	\mcF_{1/2n,\kappa/2n}(\lambda) \\
	 & = & \theta(\lambda)[G_{e}(z_{1}-z_{2}-\frac{\lambda+\kappa}{n}, 
\lambda) 
	 \tilde{\mcF}_{-\lambda/n}P - 
G_{e}(z_{1}-z_{2}-\frac{\lambda+\kappa}{n}, \kappa) 
\tilde{\mcF}_{-\kappa/n}] \\
     & & \\
	 & = & \frac{\theta'(0) \theta(z_{1}-z_{2} 
	 -\frac{\lambda+\kappa}{n}+\lambda)} 
	 {\theta(z_{1}-z_{2}-\frac{\lambda+\kappa}{n})} 
	 \tilde{\mcF}_{-\lambda/n}P - \frac{\theta(\lambda)\theta'(0) 
\theta(z_{1}-z_{2} 
	 -\frac{\lambda+\kappa}{n}+\kappa)} 
	 {\theta(z_{1}-z_{2}-\frac{\lambda+\kappa}{n}) \theta(\kappa)} 
	 \tilde{\mcF}_{-\kappa/n}.
\end{eqnarray*}  

Then $\mcR$ an entire function of $\lambda$ and $\mcR(0) = \theta'(0)P$.  
To see that $V^{e}_{n} \otimes V^{e}_{n}$ is invariant under $\mcR$, 
from the quasi-double periodicity of $\vartheta_{1}(z, \tau)$ we have:

\begin{itemize}
	\item[]  $\mcR(\lambda) \cdot f(z_{1}+1,z_{2}) = (-1)^{n-1} 
	\mcR(\lambda) \cdot f(z_{1},z_{2}) 
	= \mcR(\lambda) \cdot f(z_{1},z_{2}+1)$
	\item[]  $\mcR(\lambda) \cdot f(z_{1}+\tau,z_{2}) = 
	e^{-\pi i n \tau - 2\pi i n z_{1}} 
	\mcR(\lambda) \cdot f(z_{1},z_{2})$, 
	and 
	\item[]  $\mcR(\lambda) \cdot f(z_{1},z_{2}+\tau) = 
	e^{-\pi i n \tau - 2\pi i n z_{2}} \mcR(\lambda) \cdot f(z_{1},z_{2})$
\end{itemize}

For the holomorphicity of $\mcR(\lambda) \cdot f$ in each variable, 
the possible poles of $\mcR(\lambda) \cdot f(z_{1},z_{2})$ occur where 
$z_{1} - z_{2} -\frac{\lambda+\kappa}{n} \equiv 0$, 
modulo the lattice $\Lambda_{\tau}$.  Observing that these simple 
zeros 
of the denominator are also zeros of
$\vartheta_{1}(z_{1}-z_{2} - \frac{\lambda + \kappa}{n}) 
\mcR(\lambda) \cdot f(z_{1},z_{2})$, we see that 
$\mcR(\lambda) \cdot f(-,z_{2}), \mcR(\lambda) \cdot f(z_{1},-) \in 
{\mathcal H}({\mathbb C})$, 
and thus, $\mcR(\lambda)$ preserves $V^{e}_{n} \otimes V^{e}_{n}$.  Henceforth, 
we denote the restriction of $\mcR(\lambda)$ to $V^{e}_{n} \otimes V^{e}_{n}$ by 
$R(\lambda)$.  

To show that $R(\lambda)$ is equivalent to Belavin's $R$-matrix \cite{B}, 
$R_{B}(\lambda)$, 
define $S,T \in {\rm End}({\mathcal M}({\mathbb C}))$ by:
\begin{equation}
	\begin{array}{lll}
	(S \cdot f)(z) & = & e^{\pi i (n-1)/n}f(z + \frac{1}{n})  \\
	(T \cdot f)(z) & = & e^{\pi i \tau/n - 2\pi i z}f(z - 
\frac{\tau}{n}).
\end{array}
	\label{eq:ST}
\end{equation}

Note that $V^{e}_{n}$ is invariant under $S$ and $T$, and on $V^{e}_{n}$, 
$S^{n} = T^{n} = I$ and $TS = \omega ST$, where $\omega = e^{2\pi i/n}$.

The action of $S \otimes S$ and $T \otimes 
T$ 
on $V^{e}_{n} \otimes V^{e}_{n}$ is given by:
$$S\otimes S \cdot f(z_{1},z_{2}) = e^{2\pi i (n-1)/n}f(z_{1} + 
\frac{1}{n},
z_{2} + \frac{1}{n}) = \omega^{-1} f(z_{1} + \frac{1}{n}, z_{2} + 
\frac{1}{n})$$ and 
$$T\otimes T \cdot f(z_{1},z_{2}) = e^{2\pi i \tau/n - 2\pi i 
(z_{1}+z_{2})}
f(z_{1} - \frac{\tau}{n},z_{2}- \frac{\tau}{n}),$$
and we obtain:
\begin{itemize}
\item[{\sf (a)}] $R(\lambda)$ is {\sf completely ${\mathbb Z}_{n}$ 
	symmetric}, i.e. 
	\begin{eqnarray*}
		(S \otimes S)^{-1}R(\lambda)(S \otimes S) & = & R(\lambda) \\
		(T \otimes T)^{-1}R(\lambda)(T \otimes T) & = & R(\lambda) 
	\end{eqnarray*}
and 
\item[{\sf (b)}] $R(\lambda)$ has the following quasi-double 
periodicity: 
\begin{itemize}
	\item[{\sf (b1)}]  $R(\lambda + 1) = - (S \otimes 
	1)^{-1}R(\lambda)(S \otimes 1),$
	\item[{\sf (b2)}]  $R(\lambda + \tau) = e^{-2\pi i (\xi + 
\lambda)} 
	(T \otimes 1)R(\lambda)(T \otimes 1)^{-1}$, with $\xi = 
	-\frac{\kappa}{n} + \frac{\tau}{2} + \frac{1}{2}.$
\end{itemize}
\end{itemize}

To compute the matrix coefficients,  
$R(\lambda)_{ij}^{kl}$, with respect to the basis ${\mathcal B}$ of 
Proposition \ref{ebasis}, $S$ and $T$ act on ${\mathcal B}$ by:
\begin{eqnarray*}
		S \cdot \psi_{a}(z) & = & \omega^{a}\psi_{a}(z),  \\
	T \cdot \psi_{a}(z) & = & \psi_{a-1}(z).
\end{eqnarray*}

\noindent Hence, $R$ is given in Belavin's representation by:
$$R(\lambda) = \sum_{\alpha \in {\mathbb Z}_{n}^{\times 2}} 
{\rm w}_{\alpha}(\lambda) I_{\alpha} \otimes I_{\alpha}^{-1}$$
where $I_{\alpha} = 
S^{\alpha_{1}}T^{\alpha_{2}}$ and ${\rm w}_{\alpha} \in {\mathcal 
H}({\mathbb C})$.  
Thus, 
$$R(\lambda)_{ij}^{kl} = \delta_{i+j,k+l (n)}\sum_{\gamma = 0}^{n-1} 
{\rm 
w}_{\gamma,i-k}(\lambda) \omega^{\gamma (k-j)}.$$

From (b1) and (b2), ${\rm w}_{\alpha}(\lambda)$ has quasi-double 
periodicity:
\begin{eqnarray*}
	{\rm w}_{\alpha}(\lambda + 1) & = & -\omega^{\alpha_{2}}{\rm 
	w}_{\alpha}(\lambda),  \\
	{\rm w}_{\alpha}(\lambda + \tau) & = & e^{-2\pi i(\xi + 
	\lambda)}\omega^{\alpha_{1}}{\rm w}_{\alpha}(\lambda)  ~,~~ 
	\xi = -\frac{\kappa}{n} + \frac{\tau}{2} + \frac{1}{2}
\end{eqnarray*}
which, together with the 
initial condition $R(0) = \theta'(0)P$, uniquely determine 
$${\rm w}_{\alpha}(\lambda) = {\rm 
w}_{\alpha_{1},\alpha_{2}}(\lambda) = \frac{\theta'(0) 
\theta_{\frac{1}{2} + \frac{\alpha_{2}}{n}, 
\frac{1}{2} - \frac{\alpha_{1}}{n}}(\lambda - 
\frac{\kappa}{n},\tau)}{n \theta_{\frac{1}{2} + \frac{\alpha_{2}}{n}, 
\frac{1}{2} - \frac{\alpha_{1}}{n}}(-\frac{\kappa}{n},\tau)}$$
and 
$$R(\lambda)_{ij}^{kl} = \delta_{i+j,k+l (n)} ~ 
\frac{\theta_{\frac{1}{2},\frac{1}{2}}'(0,n\tau) 
\theta_{\frac{i-j}{n}+ \frac{1}{2},\frac{1}{2}}(\lambda - 
\kappa,n\tau)}{\theta_{\frac{i-k}{n}+ 
\frac{1}{2},\frac{1}{2}}(-\kappa,n\tau) \theta_{\frac{k-j}{n}+ 
\frac{1}{2},
\frac{1}{2}}(\lambda,n\tau)}.$$
\end{proof}

\begin{rem}  For $n$ odd, these are the same 
representation spaces found by 
Felder and Pasquier to realize Belavin's $R$-matrices.  
Our choice of representation spaces for $n$ even 
provide a concrete realization of the 
degeneration (given in the next section) from the 
Belavin $R$-matrices to the (twisted) trigonometric Cremmer-Gervais 
$R$-matrices.  
\end{rem}

At the trigonometric level,  the twisted Shibukawa-Ueno 
$R$-operators $\mcR^{t}_{0,\kappa/2n}(\lambda)$ give rise to 
the affinized Cremmer-Gervais $R$-matrices (see Section 4.2).  
However, at the elliptic level, there is no 
such analog.

\begin{thm} There are no finite dimensional  
subspaces of ${\mathcal M}({\mathbb C}^{2})$ invariant 
under $\mcR^{e}_{0,\beta}(\lambda)$.
\end{thm}

\begin{proof}
It suffices to show that $\mcR^{e}(\lambda)$ 
has no finite dimensional invariant subspaces.  
Suppose $\rho$ is an eigenvalue of 
$\mcR^{e}(\lambda) = G(z_{1}-z_{2},\lambda)P - 
G(z_{1}-z_{2},\kappa)I$.  Then 
$$\frac{f(z_{1},z_{2})}{f(z_{2},z_{1})} = 
\frac{G(z_{1}-z_{2},\lambda)}{G(z_{1}-z_{2},\kappa) + \rho}.$$
Hence, $\displaystyle{g(z_{1},z_{2}) = 
\frac{G(z_{1}-z_{2},\lambda)}{G(z_{1}-z_{2},\kappa) + \rho}}$ 
satisfies $g(z_{2},z_{1}) = g(z_{1},z_{2})^{-1}$, or, 
\begin{equation}
	\rho^{2} + [G(z,\kappa) + G(-z,\kappa)]\rho 
+ G(z,\kappa)G(-z,\kappa) - G(z,\lambda)G(-z,\lambda) = 0.
	\label{eq:eigenval}
\end{equation}

Setting $x=0$ in the 3-term equation, 
$$\theta(y)^{2}\theta(z+w)\theta(z-w) + 
\theta(z)^{2}\theta(w+y)\theta(w-y) + 
\theta(w)^{2}\theta(y+z)\theta(y-z) =0,$$
and we obtain \\
\\
$G(z,\lambda)G(-z,\lambda) - G(z,\kappa)G(-z,\kappa)$
\begin{eqnarray*}
	 & = & 
	\frac{\theta'(0)^{2}}{\theta(z)^{2} \theta(\lambda)^{2} \theta(\kappa)^{2}} 
	\left[\theta(\kappa)^{2}\theta(z+\lambda) \theta(z-\lambda)  - 
	\theta(\lambda)^{2} \theta(z + \kappa) \theta(z - \kappa) \right]  \\
	 & = & 
	 \frac{\theta'(0)^{2}}{\theta(z)^{2} \theta(\lambda)^{2} \theta(\kappa)^{2}} 
	 \left[\theta(z)^{2}\theta(\kappa +\lambda) \theta( \kappa -\lambda) \right] \\
	 & = & G(\kappa,\lambda)G(-\kappa,\lambda).
\end{eqnarray*}
That is, the constant term of (\ref{eq:eigenval}) is $z$-independent. 

The linear term of (\ref{eq:eigenval}), $G(z,\kappa) + G(-z,\kappa) = 
\frac{\theta'(0)}{\theta(\kappa)}\left[ \frac{\theta(z + \kappa) + 
\theta(z - \kappa)}{\theta(z)}\right]$ can be seen to be 
$z$-dependent 
by observing that $\theta(z + \kappa) + \theta(z - \kappa)$ 
is not a scalar multiple of $\theta(z)$, since it does not have the 
required quasi-periodicity with respect to $\tau$.
Hence, 
$\mcR^{e}(\lambda)$ has no eigenfunctions, and therefore has no 
finite dimensional invariant subspaces. 
\end{proof}

\subsection{Trigonometric degeneration}
	
For $k = 0, \ldots, 
n-1$, consider the basis of $V^{e}_{n}$ defined by
\begin{eqnarray*}
	\tilde{\psi}_{k}(z) & = & e^{-\pi i (k - 
	\frac{n-1}{2})^{2}\tau/n}\psi_{k}(z)  \\
	 & = & \sum_{m \equiv k (n)} e^{\pi i \left[ (m - 
	 \frac{n-1}{2})^{2} - (k - \frac{n-1}{2})^{2} \right] \tau/n}e^{2\pi 
	 i (m -\frac{n-1}{2}) z}.
\end{eqnarray*}
For each $k$, set 
\begin{equation}
	\phi_{k}(z)  \stackrel{\rm def}{=} 
	\lim_{{\rm Im} \tau \to \infty} \tilde{\psi}_{k}(z) = 
e^{2\pi i (k - \frac{n-1}{2})z} 
    \label{eq:phi}
\end{equation}
and define 
\begin{equation}
	V^{t}_{n} = {\rm Span} \{ \phi_{k}(z) ~|~ k = 0, \ldots, n-1\} = {\rm 
Span} 
\{ e^{2\pi i l z} ~|~ l = -j, -j + 1, \ldots, j\}
	\label{eq:Vt}
\end{equation}
where $j = \frac{n-1}{2} \in \frac{1}{2}{\mathbb Z}$.  

In contrast to the elliptic case the complete 2-parameter family of 
twisted SU operators restrict simultaneously to this finite 
dimensional subspace.

\begin{thm}
The 2-parameter twisted trigonometric Shibukawa-Ueno operators 
$\mcR^{t}_{\alpha,\beta}(\lambda)$ restrict to $V^{t}_{n} \otimes V^{t}_{n}$, 
yielding a finite 
dimensional representation which is homogeneous with respect to the 
basis above. These representations are precisely the 2-parameter 
affinized Cremmer-Gervais operators described above (\ref{eq:CGab}).
\end{thm}

\begin{proof}  It was observed in \cite{SU} that $\mcR^{t}(\lambda)$ 
    restricts to $V \otimes V$ where $V = {\rm Span}\{e^{2 \pi i k 
    z}\}_{k=0}^{n-1}$.  Since $V^{t}_{n} = e^{-\pi i (n-1)z}V$ and 
    $\mcR^{t}(\lambda)$ commutes with the action of (multiplication by) any 
    symmetric function (in particular, $e^{-\pi i (n-1)(z_{1}+z_{2})}$), 
    we see that $\mcR^{t}(\lambda)$ restricts to 
    $V^{t}_{n} \otimes V^{t}_{n}$.  Since 
    $V^{t}_{n} \otimes V^{t}_{n}$ is also invariant 
    under $\mcF_{\alpha,\beta}(\lambda)$, it follows that $V^{t}_{n} \otimes 
    V^{t}_{n}$ is invariant under  
    $\mcR^{t}_{\alpha,\beta}(\lambda)$.
    
    For the matrix coefficients, with $\theta(z) = \sin \pi z$, 
    $$\frac{1}{2 \pi i}\mcR^{t}_{\alpha,\beta}(\lambda) = 
    \frac{\eta \zeta^{-2} w_{1} - \eta^{-1} \gamma^{2} w_{2}}{\hat{\eta} 
    (\zeta^{-2} w_{1} - \gamma^{2} w_{2})} \tilde{\mcF}_{\zeta^{2}}P - 
    \frac{q \zeta^{-2} w_{1} - q^{-1} \gamma^{2} w_{2}}{\hat{q} 
	(\zeta^{-2} w_{1} - \gamma^{2} w_{2})} \tilde{\mcF}_{\gamma^{2}}$$
    where $w_{k}=e^{2 \pi i z_{k}}$, $\gamma = e^{2\pi i \beta}$, 
    $q = e^{\pi i \kappa}$, $\eta = e^{\pi i \lambda}$, 
    $\zeta = e^{2 \pi i 
    \alpha \lambda}$, and $\tilde{\mcF}_{s} \cdot f(w_{1},w_{2}) = 
    f(s^{-1}w_{1},s w_{2})$.

    With respect to the basis $\{\phi_{a}(z_{1}) \phi_{b}(z_{2})\}$ 
defined by (\ref{eq:phi}),  
the matrix coefficients, given by $R^{t}_{\alpha,\beta}(\lambda) 
\cdot \phi_{i}(z_{1})\phi_{j}(z_{2})  = 
\sum_{k,l} 
R^{t}_{\alpha,\beta}(\lambda)_{ij}^{kl}\phi_{k}(z_{1}) 
\phi_{l}(z_{2})$, 
can be found by direct calculation, or by earlier observations 
regarding homogeneous operators.  The $R$-matrices $R^{t}(\lambda)$ 
(whose 
matrix coeffients are given in \cite{SU}) are homogeneous, and the 
operator $F_{\alpha,\beta}(\lambda)$ acts diagonally on $V^{t}_{n}$ so 
that  
$F_{\alpha,\beta}(\lambda) \cdot \phi_{a}(z_{1})\phi_{b}(z_{2}) = 
{\rm exp}[2 \pi i (a-b)(\alpha \lambda - \beta)] 
\phi_{a}(z_{1})\phi_{b}(z_{2}) = \zeta^{a-b} \gamma^{b-a} 
\phi_{a}(z_{1})\phi_{b}(z_{2})$.  Thus, we obtain  
\begin{eqnarray*}
  % \frac{\hat{q}\hat{\eta}}{2 \pi i} 
    R^{t}_{\alpha,\beta}(\lambda)_{ij}^{kl} & = & 
    \zeta^{2(i-k)} \gamma^{2(j-k)} R^{t}(\lambda)_{ij}^{kl}  \\
     & = &   \frac{2 \pi i}{\hat{q}\hat{\eta}}\zeta^{2(i-k)} 
     \gamma^{2(j-k)} \begin{cases}
			\displaystyle{\widehat{q\eta^{-1}}} &  i=j=k=l,  \\
			\displaystyle{-\hat{\eta}q^{\text{sgn}(i-j)}} & 
			i=k \neq j=l, \\
			\hat{q}\eta^{\text{sgn}(j-i)} 
			& l=i \neq k=j,  \\
			\text{sgn}(j-i)\hat{q}\hat{\eta} & 
			\min(i,j) < k < \max(i,j)\\
			&\text{and }i+j=k+l, \\
			0 & {\rm otherwise}.
		\end{cases}.
\end{eqnarray*}

When $\alpha = 0$ and $\beta = \kappa/2n$ (so $\gamma^{n} = q$), 
these are the standard affinizations of the 
Cremmer-Gervais $R$-matrices.  More generally, we have (up to scalar 
multiple) $R^{t}_{\alpha,\beta}(\lambda) = 
R_{CG(\alpha,\beta-\kappa/2n)}(\lambda)$.
\end{proof}
This yields a natural degeneration of the Belavin $R$-matrix, 
$R_{B}(\lambda) = R^{e}_{1/2n,\kappa/2n}(\lambda)$, into a 
certain kind of twisted affinized Cremmer-Gervais matrix, namely 
$R_{CG(1/2n,0)}(\lambda)$.  Set 
$R^{\theta}_{1/2n,\kappa/2n} = R^{\theta}_{\kappa}$.  
>From (\ref{eq:phi}), together with
$$\lim_{{\rm Im}~\tau \to \infty} R^{e}_{\kappa}(\lambda) 
\tilde{\psi}_{a}(z_{1})\tilde{\psi}_{b}(z_{2}) = R^{t}_{\kappa}(\lambda) 
\phi_{a}(z_{1})\phi_{b}(z_{2}),$$
we obtain

\begin{thm}
The matrices 
$R_{CG(1/2n,0)}(\lambda)$ are
trigonometric degenerations of Belavin's $R$-matrices.  The 
degeneration is given by
$$\lim_{{\rm Im} \tau \to \infty} (G \otimes 
G)^{-1}R_{B}(\lambda)(G \otimes G) = 
R_{CG(1/2n,0)}(\lambda)$$
where $G$ is defined by $G_{ab} = \delta_{ab}e^{-\pi i (a - 
\frac{n-1}{2})^{2}\tau/n}$ for $a,b = 0, \ldots, n-1$.
\end{thm}  

\begin{rem}  This result is analogous to one obtained by  
Antonov, Hasegawa, and Zabrodin \cite{AHZ}. 
\end{rem}

\subsection{Rational degeneration}
The trigonometric to rational degeneration is much more straightforward 
than the elliptic to trigonometric degeneration. We omit the details which are 
easy to verify. Essentially we find that the degeneration of the 
Shibukawa-Ueno operators induces, via an appropriate representation, 
an affinized version of the degeneration of the Cremmer-Gervais R-matrices 
into the "Jordan-Cremmer-Gervais" operators described in \cite{EH1}. 

	Let us first look at twisted rational Shibukawa-Ueno operators and their representations. In the rational case, when $\theta(z)=z$, the Shibukawa-Ueno operators have the simple form:
$$
	\mcR^{r}(\lambda) \cdot f(z_1,z_2)
= \frac{z_1-z_2+\lambda}{(z_1-z_2)\lambda}f(z_2,z_1)-\frac{z_1-z_2+\kappa}{(z_1-z_2)\kappa}f(z_1,z_2)
$$
while the twisted versions have the form
\begin{equation*}
\begin{split}
	\mcR^{r}_{\alpha,\beta}(\lambda) \cdot f(z_1,z_2)
= &\frac{z_1-z_2-2(\alpha\lambda+\beta)+\lambda}{(z_1-z_2-2(\alpha\lambda+\beta))
\lambda}f(z_2+2\alpha\lambda,z_1 -2\alpha\lambda)\\
&-\frac{z_1-z_2-2(\alpha\lambda+\beta)+\kappa}{(z_1-z_2-2(\alpha\lambda+\beta))\kappa}f(z_1-2\beta,z_2+2\beta)
\end{split}
\end{equation*}
Let $V^r_n = {\rm Span}\{z^{k} \mid  k = 0, \ldots, n-1\}$. It is easily seen 
that the above operators restrict to $V^r_n\otimes V^r_n$. Moreover 
$V^r_n \otimes V^r_n$ is invariant under $\mcF_{\alpha,\beta}(\lambda)$, so that 
the restriction to these finite dimensional spaces commutes with the twisting 
operation. Fixing $n$, denote by $R^r(\lambda)$ and $R^r_{\alpha,\beta}(\lambda)$, 
respectively, the restrictions of $\mcR^{r}(\lambda)$ and 
$\mcR^{r}_{\alpha,\beta}(\lambda)$ to $V^r_n\otimes V^r_n$.
Looking explicitly at the case $\alpha=0$, we see that
$$
R^r_{0,\beta}(\lambda)\cdot z_1^iz_2^j = \frac{1}{\lambda}z_2^iz_1^j - 
\frac{1}{\kappa} (z_1-2\beta)^i(z_2+2\beta)^j
+\frac{z_2^iz_1^j-(z_1-2\beta)^i(z_2+2\beta)^j}{(z_1-z_2-2\beta)}.
$$
Set $R^r_{\beta}=\lim_{\lambda\to \infty} R^r_{0,\beta}(\lambda)$. Then 
$R^r_{\beta}$ is the general Jordan-Cremmer-Gervais operator introduced 
in \cite{EH1} and $R^r_{0,\beta}(\lambda)$ is an affinization of this 
constant operator, in the sense \cite[page 296]{KS} that
$$
R^r_{0,\beta}(\lambda)= (1/\lambda)P + R^r_{\beta}.
$$
Thus, $R^r_{\alpha,\beta}(\lambda)$ are twists 
of these affinized Jordan-Cremmer-Gervais $R$-matrices, which we denote by 
$R_{JCG(\alpha,\beta)}(\lambda)$.  

For the rational degeneration, we introduce the period $2\tau_{1}$, 
taking $\mcR^{t}(\lambda)$ to be the trigonometric Shibukawa-Ueno 
$R$-operator with $\theta(z) = \sin(\pi z/\tau_{1})$, and 
$\phi_{k}(z, \tau_{1}) = \phi_{k}(z/\tau_{1})$.  Then 
$\mcR^r_{\alpha,\beta}(\lambda) =\lim_{\tau_{1} \to \infty}
\mcR^t_{\alpha,\beta}(\lambda)$.

As in the previous section, the vector space  
$V_n^t(\tau_{1}) = {\rm Span} \{\phi_{k}(z, \tau_{1})\}$ yields a 
representation of $\mcR^t(\lambda)$. 
Consider the basis of $V_n^t(\tau_{1})$ defined by
\begin{equation}
	\tilde{\phi}_{k}(z) = \tilde{\phi}_{k}(z,\tau_{1}) = 
	\sum_{l=0}^{k} (-1)^{k-l}
\left(\frac{\tau_{1}}{2\pi i}\right)^{k} \binom{k}{l} \phi_{l}(z,\tau_{1}). 
\label{eq:rbasis}
\end{equation}
Then for each $k$, $\lim_{\tau_{1} \to \infty} 
\tilde{\phi_{k}}(z,\tau_{1}) = z^{k}$. This yields the following degeneration of the 
2-parameter affinized Cremmer-Gervais $R$-matrices into this 2-parameter family 
of Jordan-Cremmer-Gervais $R$-matrices:

\begin{thm}
The representation of $\mcR^t_{\alpha,\beta}$ on $V^t_n\otimes V^t_n$ degenerates 
as $\tau_1 \to \infty$ into the representation of $\mcR^r_{\alpha,\beta}$ 
on $V^r_n\otimes V^r_n$. This yields a degeneration of the matrix 
$R_{CG(\alpha,\beta)}$ into $R_{JCG(\alpha,\beta)}$ given by 
$$\lim_{\tau_{1} \to \infty} (H \otimes 
H)^{-1}R_{CG(\alpha, \beta)}(\lambda)(H \otimes H) 
 = R_{JCG(\alpha,\beta)}(\lambda),$$
where $H$ is defined by
$$H_{ab} = \left(\frac{\tau_{1}}{2\pi i}\right)^{b} (-1)^{b-a} \binom{a}{b},\quad a,b = 0, \ldots, n-1.$$
\end{thm}

\end{document}